\documentclass{amsart}
\usepackage{amsthm}                 
\usepackage{amsmath}                
\usepackage{amsfonts,amssymb}       
\usepackage{graphicx}               
\usepackage{float}                  
\usepackage{rotating}               
\restylefloat{figure}               
\usepackage{mathrsfs}               
\usepackage{fancyhdr}               
\usepackage[breaklinks]{hyperref}   
\usepackage[numbers]{natbib}
\newtheoremstyle{mystyle}{.25cm}{.25cm}{\it}{}{\bfseries}{.}{.5em}{\thmname{#1}\thmnumber{ #2}\thmnote{ (#3)}}
\theoremstyle{mystyle}
\newtheorem{lem}{Lemma}
\newtheorem{prop}[lem]{Proposition}
\newtheorem{thm}[lem]{Theorem}
\newtheorem{cor}[lem]{Corollary}

\theoremstyle{definition}

\newtheorem{rem}[lem]{Remark}
\newcommand{\R}{\mathbb R}

\newcommand{\Z}{\mathbb Z}
\newcommand{\N}{\mathbb N}

\newcommand{\Diff}{\mbox{\rm Diff}}

\renewcommand{\L}{\mathcal L}
\newcommand{\A}{\mathbb A}
\newcommand{\id}{\text{\rm id}}
\newcommand{\dx}{\,\text{\rm d}x}
\renewcommand{\d}{\,\text{\rm d}}
\newcommand{\dw}{\text{\rm d}}
\newcommand{\D}{\mathcal D}
\newcommand{\g}{\mathfrak{g}}

\newcommand{\Ad}{\text{\rm Ad}}
\newcommand{\ad}{\text{\rm ad}}

\renewcommand{\S}{\mathbb S}
\newcommand{\T}{\mathbb T}
\renewcommand{\phi}{\varphi}

\newcommand{\eps}{\varepsilon}

\newcommand{\ska}[2]{\left\langle #1,#2\right\rangle}

\newcommand{\set}[2]{\left\{#1;\;#2\right\}}
\renewcommand{\i}{\text{i}}

\newcommand{\bea}{\begin{eqnarray}}
\newcommand{\eea}{\end{eqnarray}}
\newcommand{\beq}{\begin{equation}}
\newcommand{\eeq}{\end{equation}}
\renewcommand{\phi}{\varphi}

\renewcommand{\autoref}[1]{\text{Eq.}~\eqref{#1}}
\begin{document}
\title{The two-dimensional periodic $b$-equation on the diffeomorphism group of the torus}
\author{Martin Kohlmann}
\address{Peter L. Reichertz Institute for Medical Informatics, University of Braunschweig, D-38106 Braunschweig, Germany}
\email{martin.kohlmann@plri.de}
\begin{abstract}
In this paper, the two-dimensional periodic $b$-equation is discussed under geometric aspects, i.e., as a geodesic flow on the diffeomorphism group of the torus $\T=\S\times\S$. In the framework of Arnold's \cite{A66} famous approach, we achieve some well-posedness results for the $b$-equation and we perform explicit curvature computations for the 2D Camassa-Holm equation, which is obtained for $b=2$. Finally, we explain the special role of the choice $b=2$ by giving a rigorous proof that $b=2$ is the only case in which the associated geodesic flow is weakly Riemannian.
\end{abstract}
\maketitle
\tableofcontents
\section{Introduction}
In the present paper, we are concerned with the two-dimensional periodic $b$-equation \cite{HS03}
\beq\label{2Db}\frac{\partial m_i}{\partial t}=-u^j\frac{\partial m_i}{\partial x^j}-m_j\frac{\partial u^j}{\partial x^i}-(b-1)m_i\frac{\partial u^j}{\partial x^j},\quad i=1,2,\quad m=\A u;\eeq
here, $u_1$ and $u_2$ are real-valued functions of space coordinates $x_1,x_2\in\R/\Z\simeq\S$ and time $t\geq 0$ and $\A$ is a linear operator whichs maps the variable $u=(u_1,u_2)$ to the momentum $m=(m_1,m_2)$. If $u_2=m_2=0$ and $u_1$ does not depend on $x_2$, \autoref{2Db} reduces to the one-parameter family
\beq\label{beq}m_t=-m_xu-bu_xm\eeq
called the $b$-equation \cite{DHH02,EY08,HS03}. Note that \autoref{beq} includes the famous Camassa-Holm (CH) equation \cite{CH93,CE98}
\beq\label{CH}u_t+3uu_x=2u_xu_{xx}+uu_{xxx}+u_{txx}\eeq
for $b=2$ and the Degasperis-Procesi (DP) equation \cite{DHH02,DP99,ELY07}
$$u_t+4uu_x=3u_xu_{xx}+uu_{xxx}+u_{txx}$$
for $b=3$ with $m=u-u_{xx}$ in both cases. The CH equation is known to describe the motion of shallow water waves over a flat bad under the action of gravity in a certain regime; the function $u(t,x)$ models the wave's height at time $t$ and position $x\in\S$. Although the DP equation has been derived in search of integrable variants of CH, it has afterward been related to ideal fluid motion as well. For an extensive presentation of the hydrodynamical relevance and integrability issues of CH and DP, we refer the reader to \cite{CL09,I05,I07,J03}.

The two-dimensional analogs of CH and DP are the denoted as 2D-CH and 2D-DP and are obtained from \eqref{2Db} for the choice $b=2$ and $b=3$ respectively, with $\A=\text{diag}(1-\Delta,1-\Delta)$. Since the two-dimensional Camassa-Holm equation is particularly important for the paper at hand, we now shed some light on its physical meaning and the features of its solutions which exhibit typically nonlinear phenomena as wave-breaking or peakons.

First, 2D-CH is related to ideal fluid motion \cite{KSD01} and it can be shown that the terms on the right-hand side model convection, stretching and expansion of a fluid with velocity $u=(u_1,u_2)$ and momentum $m=(m_1,m_2)$, cf.~\cite{HM05,HSS09},
\beq\label{2DCH}m_t=-\underbrace{u\cdot\nabla m}_{\text{\rm convection}}-\underbrace{(\nabla u)^Tm}_{\text{\rm stretching}}-\underbrace{(\nabla\cdot u)m}_{\text{\rm expansion}},
\quad
m=
\left(
\begin{array}{cc}
1-\Delta & 0 \\
0 & 1-\Delta \\
\end{array}
\right)
u.
\eeq
The 2D-CH equation is also known as the Navier-Stokes-$\alpha$ model \cite{BLY09,V02}, since its viscous variant is closely related to the two-dimensional Navier-Stokes system. Furthermore, the viscous Camassa-Holm equation is an appropriate model for turbulent channel and pipe flow; see, e.g., \cite{CFHOTW98,CFHOTW99,CFHOTW99b,FHT02} where the authors present two- and three-dimensional variations of \autoref{CH}. A derivation of higher-dimensional Camassa-Holm equations from the physical point of view can also be found in \cite{H99,HMR98b} and integrability properties are discussed in \cite{KZ99}. It is important to remark that \autoref{2DCH} also possesses an interesting specialization that takes it into a modified two-component Camassa-Holm system (MCH2) \cite{HNT09}
\bea
\begin{array}{rcl}
q_t+vq_x+2qv_x+\rho(1-\partial_x^2)^{-1}\rho_x & = & 0,\nonumber\\
\rho_t+(\rho v)_x & = & 0,\nonumber
\end{array}
\eea
where $q$, $v$ and $\rho$ depend on time $t$ and the single space variable $x$. By setting $x_1=x$, $u=(v,(1-\partial_x^2)^{-1}\rho)$ and $m=(q,\rho)$ in \autoref{2DCH}, we obtain the MCH2 system from the 2D-CH equation; see also \cite{HI11}.

Second, the mathematical theory behind \autoref{2DCH} is rich and interesting \cite{HMR98,HM05,HSS09}: The 2D-CH can be derived by means of a Lagrangian formulation. Let $\ell(u)=\frac{1}{2}\int_\T u\cdot m\d^2x$ be the kinetic energy Lagrangian. The 2D-CH equation is the Euler-Lagrange equation obtained from the variational principle
$$\delta\int\ell(u)\d t=0$$
using variations of the form $\delta u=\dot v+v\cdot\nabla u-u\cdot\nabla v$. Via the Legendre transformation, the Hamiltonian formulation of the 2D-CH reads
$$m_t=\{m,H\},\quad H(m)=\ska{m}{u}-\ell(u),$$
where the duality pairing $\ska{\cdot}{\cdot}$ of momentum and velocity is given by integration and $\{\cdot,\cdot\}$ is the so-called ideal fluid Lie-Poisson bracket \cite{MR99}. Next, the 2D-CH can be nicely rewritten in the form
$$m_t+\L_um=0,$$
where $\L$ denotes the Lie derivative of the momentum one form density $m=m_i\dx^i\otimes\d^2x$ with respect to the velocity vector field $u=u^i\partial_{x^i}$. Its $N$-peakon solutions are given by the measure-valued (that is, density valued) ansatz
$$m(t,x)=\sum_{a=1}^N\int P^a(t,s)\delta(x-Q^a(t,s))\d s.$$
Similarly to the one-dimensional Camassa-Holm equation \cite{CK02,CK03,K99,S98}, it is known that the 2D version \eqref{2DCH} is related to a geodesic flow on the diffeomorphism group of the torus $\T=\S\times\S$, i.e., \autoref{2DCH} belongs to the so-called EPDiff equations \cite{GB09,GB11}. It goes back to Arnold's fundamental observation \cite{A66} and the famous work \cite{EM70} of Ebin and Marsden that the geometric viewpoint on the one-dimensional CH equation is similar to the geometric picture of the rotational motion of a rigid body in $\R^3$: In both cases the dynamics is equivalent to geodesic motion on an appropriate Lie group with respect to a certain invariant metric or a compatible affine connection respectively. The 2D Camassa-Holm equation \eqref{2DCH} can also be captured within this powerful approach which is explained more detailed in Appendix~\ref{app_geom2} of the present paper. The geometric theory is not only aesthetically appealing but has also important applications concerning well-posedness and stability issues; see, e.g., \cite{GB09,GB11} where the author proves well-posedness of the $n$-dimensional Camassa-Holm equation on compact Riemannian manifolds with certain boundary conditions, and \cite{EKL11,G00,GO06,M98,M02} for an extensive presentation of related material.

Our paper contains the following new results: We show that \autoref{2Db} re-expresses geodesic motion on the diffeomorphism group of $\T$ for \emph{any} choice of the parameter $b$. Our approach will be carried out on a scale of Sobolev spaces whose regularity properties are presented very detailed in a novel paper of Inci, Kappeler and Topalov \cite{IKT11}. For $b=2$, we will also be able to study geodesics on the Lie group of smooth diffeomorphisms of $\T$, associated with \autoref{2DCH}. In this context, we derive well-posedness in the Sobolev spaces $H^s(\T)$ for \autoref{2Db} and a general $b$ and for \autoref{2DCH} in the smooth category.
\thm[Well-posedness in $H^s(\T)$]\label{thm_lwp1} Let $s>3$ and $b\in\R$. There is an open neighborhood $U\subset H^s(\T)$ of zero, such that for any $u_0\in U$ there is $T>0$ and a unique solution
$$u\in C([0,T);H^s(\T))\cap C^1([0,T); H^{s-1}(\T))$$
to the two-dimensional $b$-equation \eqref{2Db} satisfying the initial condition $u(0)=u_0$ and depending continuously on $u_0$, i.e., the mapping
$$u_0\mapsto u,\quad U\to C([0,T);H^s(\T))\cap C^1([0,T); H^{s-1}(\T))$$
is continuous.
\endthm\rm
\thm[Well-posedness in $C^{\infty}(\T)$]\label{thm_lwp2} Let $b=2$. There exists an open neighborhood $U\subset H^4(\T)$ of zero such that for any $u_0\in U\cap C^{\infty}(\T)$ there exists $T>0$ and a unique solution $u\in C^{\infty}([0,T),C^{\infty}(\T))$ of \autoref{2Db} satisfying the initial condition $u(0)=u_0$ and depending smoothly on $u_0$, i.e., the map
$$u_0\mapsto u,\quad U\cap C^{\infty}(\T)\to C^{\infty}([0,T);C^{\infty}(\T))$$
is smooth.
\endthm\rm
Next, we aim to give a proof of the remarkable property that the geodesic flow associated with \autoref{2Db} is (weakly) Riemannian, with respect to a special class of inertia operators, if and only if $b=2$. Since the \emph{if}-part has already been established in, e.g., \cite{GB09,GB11}, and is also recalled in Appendix~\ref{app_geom1} for convenience, we complement the \emph{only if}-part in this paper. In particular, we obtain that there is no right-invariant metric belonging to the class under consideration which is associated with the 2D-DP equation on the diffeomorphism group of the torus. Note that a similar result for the one-dimensional $b$-equation has been established in \cite{ES10,K09}. Our main theorem in this context is
\thm\label{thm_b=2} Let $b\geq 2$ be an integer and $\mathbb L=\text{\rm diag}(1-\Delta,1-\Delta)$. Suppose that there is a regular
inertia operator $\A=\text{\rm diag}(A,A)$, $A\in\mathcal L_{\text{\rm is}}^{\text{\rm sym}}(C^{\infty}(\T)),$ such that the 2D $b$-equation
$$m_t=-u\cdot\nabla m-(\nabla u)^Tm-(b-1)m(\nabla\cdot u),\quad m=\mathbb L u,$$
is the Euler equation on $\Diff^{\infty}(\T)$ with respect to the right-invariant metric $\rho_\A$ defined by $\A$. Then $b=2$ and $\A=\mathbb L$.\endthm
\endthm\rm
With the right-invariant metric on the diffeomorphisms $\T\to\T$ for $b=2$, we can define a sectional curvature $S$ associated with \autoref{2DCH} and we derive a convenient formula showing that $S$ can be expressed in terms of the Christoffel operator for \autoref{2DCH}. As an application, we write down a large number of two-dimensional spaces on which $S$ is positive. These results are inspired by the work done for other variants of the Camassa-Holm equation in \cite{EKL11,M98,M02}.

\begin{thm}\label{thm_curv} Let $R$ be the curvature tensor on
the torus diffeomorphism group equipped with the right-invariant metric $\ska{\cdot}{\cdot}$ associated with the 2D Camassa-Holm equation. Then the sectional curvature $S(u,v) :=\ska{R(u,v)v}{u}$ is given at the identity by
\begin{equation}\label{Suvexpression}
  S(u,v)=\ska{\Gamma(u,v)}{\Gamma(u,v)}-\ska{\Gamma(u,u)}{\Gamma(v,v)}+R(u,v),
\end{equation}
where
\bea R(u,v)\!\!\!&=&\!\!\!\ska{\nabla u\cdot u}{\nabla v\cdot v}-\ska{\nabla u\cdot v}{\nabla u\cdot v}+\ska{\nabla v\cdot u}{\nabla u\cdot v}-\ska{\nabla v\cdot u}{\nabla v\cdot u}\nonumber\\
&&+\ska{[\nabla(\nabla u\cdot u)]\cdot v}{v}-\ska{[\nabla(\nabla u\cdot v)]\cdot v}{u}+\ska{[\nabla(\nabla v\cdot u)]\cdot v}{u}\nonumber\\
&&-\ska{[\nabla(\nabla v\cdot u)]\cdot u}{v}-\ska{\nabla v(\nabla u\cdot u)}{v}-\ska{\nabla u(\nabla v\cdot v)}{u}\nonumber\\
&&+\ska{\nabla v(\nabla v\cdot u)}{u}+\ska{\nabla u(\nabla v\cdot u)}{v}.\nonumber\eea
In particular, $S$ is positive on any of the two-dimensional spaces spanned by the canonical basis vectors $e_1$, $e_2$ and
$$
v=
\left(
    \begin{array}{c}
      \sin (k_1x)\sin (k_2y) \\
      \sin (k_1x)\sin (k_2y)\\
    \end{array}
\right),\quad k_1,k_2\in\N.
$$
\end{thm}

Our paper is organized as follows: In Section~\ref{sec_preliminaries} we first recall some elementary facts about the diffeomorphism group of $\T$ and some results of \cite{IKT11}. In Section~\ref{sec_2DCH} we develop the geometric picture for the two-dimensional $b$-equation in the $H^s$-category and proceed with a proof of the results in Theorem~\ref{thm_lwp1} and Theorem~\ref{thm_lwp2} using the geometric theory. Next, we compute the curvature of the diffeomorphism group associated with 2D-CH and give a proof of Theorem~\ref{thm_curv}. Section~\ref{sec_b} is devoted to the special case $b=2$ of the two-dimensional $b$-equation; here we work out the proof of Theorem~\ref{thm_b=2}. Finally, in Section~\ref{sec_outlook}, we discuss some open problems and further tasks; in particular, we present the so-called $\mu$-variant of the periodic two-dimensional $b$-equation which has not been discussed in the literature up to now. There are two appendices presenting some background information concerning the geometric picture for the 2D-CH equation and the commonalities in the geometric formalism for the rigid body, the one-dimensional CH equation and the two-dimensional variant \autoref{2DCH}.
\\[.25cm]\emph{Acknowledgement.} The author thanks Helge Gl\"ockner (University of Paderborn) and his group for asking about the issues discussed in the present paper. A cordial thank goes to the anonymous referees for several helpful remarks that improved the preliminary version of the paper.
\section{The diffeomorphism group of the torus: Some preliminary facts}\label{sec_preliminaries}
Let $\S=\R/\Z$ and denote by $\T=\S^2$ the two-dimensional torus, i.e., $\T$ consists of equivalence classes of pairs of real numbers such that $n=(n_1,n_2)$ and $m=(m_1,m_2)$ are equivalent if and only if $n_i-m_i$ is an integer, for $i\in\{1,2\}$. It is well known that $\T$ is a smooth manifold with trivial tangent bundle $T\T\simeq\T\times\R^2$. We will label the coordinates in $\T$ by $x=x_1$ and $y=x_2$. For any differentiable function $u=(u_1,u_2)$ on $\T$ we write $\nabla u=(\nabla u_1,\nabla u_2)$, where $\nabla=(\partial_x,\partial_y)$ is the nabla operator. Sobolev mappings between two smooth manifolds $M$ of order $n$ and $N$ of order $d$, for $s>n/2$, are defined as follows: A continuous mapping $f\colon M\to N$ is an element of $H^s(M,N)$ if for any two charts $\psi\colon\mathcal U\to U\subset\R^n$ and $\eta\colon\mathcal V\to V\subset\R^d$ with $f(\mathcal U)\subseteq\mathcal V$, $\eta\circ f\circ\psi^{-1}\colon U\to V$ is an element of the usual Sobolev space $H^s(U,\R^d)$. Let $H^s(\T)=H^s(\T,\R^2)$. Furthermore, we introduce the diffeomorphism groups
$$\D^s(\T)=\set{\phi\in\Diff^1_+(\T)}{\phi\in H^s(\T)},\quad s>2,$$
where $\Diff_+^1(\T)$ denotes the orientation-preserving $C^1$-diffeomorphisms of $\T$; precisely, $\phi\in\Diff_+^1(\T)$ if and only if $\phi\in C^1(\T;\T)$ is a diffeomorphism and $|\nabla\phi|=\det\nabla\phi>0$. Elements of $\D^s(\T)$ are referred to as orientation-preserving $H^s$ diffeomorphisms. It is shown in \cite{IKT11} that $\D^s(\T)$ is open in $H^s(\T)$ and hence is a $C^{\infty}$-Hilbert manifold. Moreover, $\D^s(\T)$ is a topological group, i.e., the group product $(\phi,\psi)\mapsto\phi\circ\psi$ and the inversion map $\phi\mapsto\phi^{-1}$ are continuous maps $\D^s(\T)\times\D^s(\T)\to\D^s(\T)$ and $\D^s(\T)\to\D^s(\T)$ respectively. For $s\to\infty$, the groups $\D^s(\T)$ approximate the Lie group
$$\Diff^{\infty}(\T)=\set{\phi\in\Diff^1_+(\T)}{\phi\in C^{\infty}(\T;\T)}=\bigcap_{s>2}\D^s(\T),$$
a $C^{\infty}$-Fr\'echet manifold on which multiplication and inversion are smooth maps. Let $\phi\in\D^s(\T)$ and let $\gamma(t)\subset\D^s(\T)$ be a $C^1$-curve starting at $\phi$. Then the vector $\dot\gamma(0)$ is an element of $H^s(\T)$ and hence
$$T_\phi\D^s(\T)\simeq H^s(\T)\simeq\set{u\cdot\nabla}{u\in H^s(\T)};$$
the latter space denotes the $H^s$ vector fields on the torus. We further have the identification
$$T\D^s(\T)=\bigcup_{\phi\in\D^s(\T)}\{\phi\}\times T_\phi\D^s(\T)\simeq\D^s(\T)\times H^s(\T)$$
and $T_\phi\D^s(\T)^*\simeq H^s(\T)$ for any $\phi\in\D^s(\T)$. Let $I_\phi\colon\D^s(\T)\to\D^s(\T)$ denote the inner automorphism $\psi\mapsto\phi\circ\psi\circ\phi^{-1}$. We let
$$\Ad_\phi v=(D_\id I_\phi)v=[(\nabla\phi)\cdot v]\circ\phi^{-1},\quad\forall(\phi,v)\in\D^s(\T)\times H^s(\T),$$
and
$$\ad_uv=\left.\frac{\d}{\d t}\right|_{t=0}\Ad_{\phi(t)}v=\nabla u\cdot v-\nabla v\cdot u,$$
where $\phi(t)$ is a $C^1$-curve in $\D^s(\T)$ starting at $\id$ with $\dot\phi(0)=u$ and $v\in H^s(\T)$. The map $\ad$ satisfies the Lie bracket properties and we write
\beq\label{commutator}[u,v]=\nabla u\cdot v-\nabla v\cdot u\eeq
in the following. Let $\A$ be some topological isomorphism of $H^s(\T)$ such that
$$\ska{u}{v}=\int_\T u\A v\d(x,y)$$
is a scalar product on $H^s(\T)$. Some easy computations show that the adjoint operators $\Ad^*_\phi$ (for the $L_2$-pairing) and $\ad_u^*$ (with respect to the scalar product $\ska{\cdot}{\cdot}$) are given explicitly by
$$\Ad_\phi^*w=(\nabla\phi)^T(w\circ\phi)|\nabla\phi|,\quad w\in H^s(\T),$$
and
\beq\ad_u^*w=\A^{-1}\left\{(\nabla u)^T\cdot\A w+\nabla(\A w)\cdot u+(\nabla\cdot u)\A w\right\},\quad w\in H^s(\T).\label{ad}\eeq
Finally, let $R_\psi,L_\psi\colon\D^s(\T)\to\D^s(\T)$, $R_\psi\colon\eta\to\eta\circ\psi$, $L_\psi\colon\eta\to\psi\circ\eta$ denote the right and left translation map on $\D^s(\T)$ respectively. It is easy to derive that
\beq\label{invariance}(D_\phi R_{{\phi}^{-1}})v=v\circ\phi^{-1},\quad (D_\phi L_{{\phi}^{-1}})v=(\nabla\phi)^{-1}v,\quad\forall v\in T_\phi\D^s(\T)\simeq H^s(\T).\eeq
The formulas \eqref{invariance} will be of great importance for the issues of this paper since we will consider right-invariant vector fields and metrics on the group $\D^s(\T)$. In case of right invariance, the sign of the commutator bracket \eqref{commutator} and the operators depending linearly on it will have to be changed.
\section{Geometry and the solutions of the two-dimensional $b$-equation}\label{sec_2DCH}
In this section, we develop the geometric picture for \autoref{2Db} by the line of arguments in \cite{HM05} and apply the geometric formalism to obtain well-posedness results for the two-dimensional periodic $b$-equation and explicit curvature computations.
\subsection{The geometric picture}
Let $\A=\text{diag}(1-\Delta,1-\Delta)$. It follows from a standard Fourier representation argument that $\A$ is a topological isomorphism $H^s(\T)\to H^{s-2}(\T)$. We can rewrite the two-dimensional $b$-equation \eqref{2Db} as
\beq\label{2DCHEuler} u_t=-\A^{-1}\{u\cdot\nabla(\A u)+(\nabla u)^T\A u+(b-1)\A u(\nabla\cdot u)\}.\eeq
Introducing the quadratic operators
$$B(u,v)=-\A^{-1}\{v\cdot\nabla(\A u)+(\nabla v)^T\A u+(b-1)\A u(\nabla\cdot v)\}$$
and
\beq\Gamma(u,v)=\frac{1}{2}\left(\nabla u\cdot v+\nabla v\cdot u+B(u,v)+B(v,u)\right),\label{Christoffel}\eeq
we find that
$$u_t+\nabla u\cdot u=\Gamma (u,u).$$
Note that $\Gamma\colon H^s(\T)\times H^s(\T)\to H^s(\T)$ which follows from the identity
$$u\cdot\nabla(\A v)-\A(\nabla v\cdot u)=\nabla v\cdot\Delta u+2\nabla v_x\cdot u_x+2\nabla v_y\cdot u_y.$$
For $u\in H^s(\T)$, let $\phi$ be the solution of the initial value problem
\bea
\left\{
\begin{array}{ccll}
\phi_t(t,z)  &=& u(t,\phi(t,z)), & z\in\T, t>0, \\
\phi(0,z)    &=& z, & z\in\T
\end{array}
\right.
\nonumber\eea
in $\D^s(\T)$ on some time interval $[0,T)$. We will use the short hand notation $\phi_t=u\circ\phi$ to emphasize that $\phi_t$ is the right-invariant vector field on $\D^s(\T)$ with value $u$ at $\id$; equivalently, $\phi$ is a local flow for the vector field $u\cdot\nabla$ on $\T$. Next, we introduce the mapping
$$\Gamma_\phi\colon\D^s(\T)\times H^s(\T)\times H^s(\T)\to H^s(\T),\quad (\phi,U,V)\mapsto\Gamma(U\circ\phi^{-1},V\circ\phi^{-1})\circ\phi$$
which is the right-invariant extension of the map $\Gamma$ on the group $\D^s(\T)$. We will call the map $\Gamma_\phi$ the \emph{Christoffel operator} for the two-dimensional $b$-equation. Moreover, the map $\phi\mapsto\Gamma_\phi$ is smooth which can be shown as explained in Appendix~\ref{app_geom1}. It follows that \autoref{2DCHEuler} is equivalent to
\beq\phi_{tt}=\Gamma_\phi(\phi_t,\phi_t)\label{2DCHgeo}\eeq
and it is known that \autoref{2DCHgeo} is a geodesic equation for the Lagrangian variable $\phi\in\D^s(\T)$ corresponding to the smooth affine connection $\bar\nabla$ on $\D^s(\T)$ given by
\beq\label{connection}\bar\nabla_XY=DY\cdot X-\Gamma(X,Y),\eeq
for any pair of smooth vector fields $X,Y$ on $\D^s(\T)$, \cite{L95}. We will call the Lagrangian variable $\phi$ the \emph{geodesic flow} for the two-dimensional $b$-equation in the following and \autoref{2DCHgeo} is referred to as the associated geodesic equation in local coordinates.

For $b=2$, the geometric picture is slightly more involved: Here, the connection $\bar\nabla$ is Riemannian in the sense that it is compatible with the metric
$$\ska{u}{v}=\int_\T u\A v\d(x,y),\quad u,v\in H^s(\T)$$
which is extended by right-invariance
\beq\label{metric}\ska{X}{Y}_\phi=\ska{X(\phi)\circ\phi^{-1}}{Y(\phi)\circ\phi^{-1}},\quad\phi\in\D^s(\T),\eeq
for any pair $X,Y$ of smooth vector fields on $\D^s(\T)$. Moreover, we have that $B(u,v)=-\ad^*_vu$ with $\ad^*$ as in \eqref{ad}. In \cite{GB09}, we find the following statements. Some more details and a proof using the notation of the present paper are provided for the convenience of the reader in Appendix~\ref{app_geom1}.
\prop\label{prop_geometry} Fix $b=2$. Let $\D^s(\T)$ for $s>3$ be the diffeomorphism group of the torus, $\mathfrak{X}(\D^s(\T))$ the space of smooth vector fields on $\D^s(\T)$ and let $\bar\nabla$ and $\ska{\cdot}{\cdot}$ be as in \eqref{connection} and \eqref{metric}. Then:
\begin{enumerate}
\item $\bar\nabla\colon\mathfrak{X}(\D^s(\T))\times\mathfrak{X}(\D^s(\T))\to\mathfrak{X}(\D^s(\T))$ is a smooth torsion-free affine connection on $\D^s(\T)$, i.e.,
\begin{itemize}
\item[(i)] $\bar\nabla_{fX+gY}Z=f\bar\nabla_XZ+g\bar\nabla_YZ$,
\item[(ii)] $\bar\nabla_X(Y+Z)=\bar\nabla_XY+\bar\nabla_XZ$,
\item[(iii)] $\bar\nabla_X(fY)=f\bar\nabla_XY+X(f)Y$,
\item[(iv)] $\bar\nabla_XY-\bar\nabla_YX=[X,Y]$,
\end{itemize}
for all $X,Y,Z\in\mathfrak X(D^s(\T))$ and all $f,g\in C^{\infty}(\D^s(\T);\R)$, and the map $\phi\mapsto(\phi,(\bar\nabla_XY)(\phi))$, $\D^s(\T)\to T\D^s(\T)$ is smooth for any $X,Y\in\mathfrak{X}(\D^s(\T))$.
\item The right-invariant inner product $\ska{\cdot}{\cdot}\colon\mathfrak{X}(\D^s(\T))\times\mathfrak{X}(\D^s(\T))\to\R$ is a weak Riemannian metric on $\D^s(\T)$.
\item The connection $\bar\nabla$ and the metric $\ska{\cdot}{\cdot}$ are compatible in the sense that
$$X\ska{Y}{Z}=\ska{\bar\nabla_XY}{Z}+\ska{\bar\nabla_XZ}{Y},\quad\forall X,Y,Z\in\mathfrak{X}(\D^s(\T)).$$
\end{enumerate}
\endprop
\rem The metric $\ska{\cdot}{\cdot}$ is only a weak Riemannian metric in the sense that the natural topology on any $T_\phi\D^s(\T)\simeq H^s(\T)$ is stronger than the topology induced by $\ska{\cdot}{\cdot}_\phi$; see \cite{EM70} for more details.
\endrem
\subsection{Well-posedness}
Since the periodic two-dimensional $b$-equation has a smooth geodesic spray, we immediately obtain well-posedness for the geodesic equation \eqref{2DCHgeo} in $H^s$.
\prop\label{prop_LWPgeo1} Let $\Gamma$ be the Christoffel map defined in \eqref{Christoffel}. There is an open neighborhood $U\subset H^s(\T)$ containing zero, for $s>3$, such that for any $u_0\in U$ the initial value problem
\bea
\left\{
\begin{array}{rcl}
\phi_{tt} & = & \Gamma_\phi(\phi_t,\phi_t), \\
\phi_t(0) & = & u_0,\\
\phi(0)   & = & \id
\end{array}
\right.
\eea
for the geodesic flow corresponding to the two-dimensional $b$-equation on $\D^s(\T)$ has a unique solution $\phi\in C^{\infty}([0,T);\D^s(\T))$, for some $T>0$, which depends smoothly on time and the initial value, i.e., the local flow $(t,u_0)\mapsto\phi$, $[0,T)\times U\to\D^s(\T)$ is smooth.
\endprop\rm
Since $\D^s(\T)$ is a topological group, we may set $u=\phi_t\circ\phi^{-1}$ to obtain a solution of \eqref{2Db} with the regularity properties specified in Theorem~\ref{thm_lwp1}. To obtain a proof of Theorem~\ref{thm_lwp2}, we show that the geodesic flow preserves its spatial regularity as we increase the regularity of the initial data, for $b=2$.
\prop Let $b=2$ and let $\Gamma$ be the Christoffel map defined in \eqref{Christoffel}. There is an open neighborhood $U_4\subset H^4(\T)$ containing zero such that for any $u_0\in U_4\cap C^{\infty}(\T)$, the initial value problem
\bea
\left\{
\begin{array}{rcl}
\phi_{tt} & = & \Gamma_\phi(\phi_t,\phi_t), \\
\phi_t(0) & = & u_0,\\
\phi(0)   & = & \id
\end{array}
\right.
\eea
for the geodesic flow corresponding to the two-dimensional Camassa-Holm equation has a unique solution $\phi\in C^\infty([0,T);\Diff^\infty(\T))$, for some $T>0$, which depends smoothly on time and the initial value, i.e., the local flow $(t,u_0)\mapsto\phi$, $[0,T)\times U_4\cap C^{\infty}(\T)\to\Diff^\infty(\T)$ is smooth.
\endprop\rm
\proof Let $\Phi\colon[0,T_4)\times U_4\to\D^4(\T)$, $\Phi(t,u_0)=\phi$, be the local flow for the 2D-CH equation on $\D^4(\T)$, obtained in Proposition~\ref{prop_LWPgeo1}. Pick $u_0\in U_4\cap C^{\infty}(\T)$ and let $[0,T_s)$ denote the interval of existence for the associated flow $\Phi(t,u_0)$ in $\D^s(\T)$, for $s\geq 4$. By uniqueness, we know that $T_s\leq T_4$ and it follows from Theorem 12.1 in \cite{EM70} that $T_s<T_4$ is not possible. Finally, an application of Lemma 3.10 in \cite{EKK11} achieves the proof.
\endproof
Theorem~\ref{thm_lwp2} is an immediate consequence of the above proposition and the fact that the torus diffeomorphism group is a Lie group in the smooth category.
\subsection{The curvature tensor of $\D^s(\T)$ associated with the 2D Camassa-Holm equation}
The existence of a smooth connection
$\bar\nabla$ on a Banach manifold $M$ immediately implies the existence of
a smooth curvature tensor $R$ defined by
$$R(X,Y)Z=\bar\nabla_X\bar\nabla_YZ-\bar\nabla_Y\bar\nabla_XZ-\bar\nabla_{[X,Y]}Z,$$
where $X,Y,Z$ are smooth vector fields on $M$, cf.~\cite{L95}. In
the case of \autoref{2DCH}, since there exists a metric $\ska{\cdot}{\cdot}$, we can also define an (unnormalized) sectional curvature
$S$ by
$$S(X,Y) := \langle R(X,Y)Y, X \rangle.$$
In this section, we will derive a convenient formula for $S$ and
use it to determine large subspaces of positive curvature for the
2-dimensional Camassa-Holm equation.

Curvature computations have a long tradition in the geometric theory of partial differential equations, see, e.g., \cite{F88,McK82}, and two-dimensional subspaces on which the sectional curvature $S$ is positive are of particular interest since the positivity of $S$ is related to stability properties of the geodesic flow, cf.~\cite{A89}. Large subspaces of positive curvature for the one-dimensional Camassa-Holm equation and its supersymmetric extension have been found in \cite{EKL11,M98,M02} by evaluating $S$ on trigonometric functions. Here we extend this discussion for the 2D Camassa-Holm equation.\\[.25cm]
\emph{Proof of Theorem~\ref{thm_curv}}. Let $U, V, W \in T_{\phi}\Diff^{\infty}(\T)$ be three tangent vectors at a point $\phi \in \Diff^{\infty}(\T)$.
The curvature tensor $R$ is given locally by
\begin{align*}
R_\phi(U, V)W = &\; D_1\Gamma_{\phi}(W, U)V - D_1\Gamma_{\phi}(W, V)U
         \\
 &  + \Gamma_{\phi}(\Gamma_{\phi}(W, V), U) - \Gamma_{\phi}(\Gamma_{\phi}(W, U), V)
\end{align*}
where $\Gamma$ is the 2D-CH Christoffel map defined in
(\ref{Christoffel}) and $D_1$ denotes differentiation with respect to
$\phi$:
$$D_1\Gamma_{\phi}(W, U)V = \frac{\dw}{\dw\eps}\bigg|_{\eps = 0} \Gamma_{\phi + \eps V}(W, U).$$
By right invariance,
$$R_\phi(X,Y)Z\circ\phi^{-1}=R_{\id}(u,v)w,$$
where $u=X(\phi)\circ\phi^{-1}$, $v=Y(\phi)\circ\phi^{-1}$ and $w=Z(\phi)\circ\phi^{-1}$, it suffices to discuss the curvature tensor at the identity. To simplify our notation, we will omit the index $\id$ of the geometric objects under discussion in the following.
First, we note that
$$D_1\Gamma(w,u)v=-\Gamma(\nabla w\cdot v,u)-\Gamma(\nabla u\cdot v,w)+\nabla\Gamma(w,u)\cdot v.$$
Thus,
\begin{align}\nonumber
S(u,v) =&\; \ska{\Gamma(\Gamma(v,v),u)}{u}-\ska{\Gamma(\Gamma(v,u),v)}{u}
    \\ \label{Suv}
& + \langle \nabla\Gamma(v,u)\cdot v, u \rangle - \langle \nabla\Gamma(v,v)\cdot u, u \rangle
    \\ \nonumber
&+\ska{-\Gamma(\nabla v\cdot v,u)-\Gamma(\nabla u\cdot v,v)+2\Gamma(\nabla v\cdot u,v)}{u}.
\end{align}
A lengthy but tedious computation similar to the calculation in the proof of Proposition 5.1 in \cite{EKL11} shows that \autoref{Suv} becomes
\begin{align} \nonumber
S(u,v)=&\;\ska{\Gamma(u,v)}{\Gamma(u,v)}-\ska{\Gamma(u,u)}{\Gamma(v,v)}
    \\ \label{Suv2}
&+\ska{\nabla u\cdot u}{\Gamma(v,v)}-\ska{
\nabla u\cdot v}{\Gamma(v,u)}
    \\ \nonumber
&+\ska{-\Gamma(\nabla v\cdot v,u)-\Gamma(\nabla u\cdot v,v)+2\Gamma(\nabla v\cdot u,v)}{u}
\end{align}
and by the definition of $\Gamma$, we obtain (\ref{Suvexpression}) after a further lengthy calculation. This proves the first assertion of Theorem~\ref{thm_curv}.
Let $\{e_1,e_2\}$ be the canonical basis of $\R^2$, let $k_1,k_2\in 2\pi\N$ and let
$$
v=
\left(
    \begin{array}{c}
      \sin (k_1x)\sin (k_2y) \\
      \sin (k_1x)\sin (k_2y)\\
    \end{array}
\right).
$$
In the following, we will make use of the identity
$$(1-\Delta)^{-1}\sin(\alpha x)\cos(\beta y)=\frac{\sin(\alpha x)\cos(\beta y)}{1+\alpha^2+\beta^2}$$
and the trigonometric formulas
$$\int_\S\cos^2(\alpha x)\dx=\int_\S\sin^2(\alpha x)\dx=\frac{1}{2},\quad\int_\S\cos(\alpha x)\sin(\beta x)\dx=0,\quad\forall\alpha,\beta\in 2\pi\N.$$
We claim that $S(e_i,v)>0$ for $i=1,2$. First we observe that for general $w\in C^{\infty}(\T)$, we have that $R(e_i,w)=0$, $i=1,2$; this follows from a straightforward computation using integration by parts. Hence
$$S(e_i,v)=\ska{\Gamma(e_i,v)}{\Gamma(e_i,v)}.$$
We leave it to the reader to perform explicitly the easy calculations leading to the formulas
\bea
S(e_1,v)&=&\frac{1}{8}\frac{2k_1^2+k_2^2}{1+k_1^2+k_2^2},\nonumber\\
S(e_2,v)&=&\frac{1}{8}\frac{2k_2^2+k_1^2}{1+k_1^2+k_2^2}.\nonumber
\eea
It follows that, for any $k_1,k_2\in 2\pi\N$, $S>0$ on the spaces $\text{span}\{e_i,v\}$, $i=1,2$. Altogether, we have completed the proof of Theorem~\ref{thm_curv}.
\hfill$\square$
\rem The curvature formula presented in \cite{EKL11} is formally identical to \autoref{Suvexpression} with $R\equiv 0$; simply the metric and the Christoffel operator have to be exchanged by their 1D two-component analogs. It is not known to the author whether $R\equiv 0$ in \eqref{Suvexpression}, but it seems to be impossible to handle the large number of terms in a calculation by hand; it is rather worthwhile using computer software to simplify the terms in $R(u,v)$, but this will not be part of this analytical paper. Observe that identities of type $S(u,v)=\ska{\Gamma(u,v)}{\Gamma(u,v)}-\ska{\Gamma(u,u)}{\Gamma(v,v)}$ do \emph{not} follow from the general theory of right-invariant metrics on Banach manifolds; the example of a modified Hunter-Saxton equation in \cite{KLM08} shows that the unnormalized sectional curvature $S$ for this equation contains an additional nontrivial term.\endrem\rm
\section{The special role of the case $b=2$}\label{sec_b}
We have shown that the 2D $b$-equation re-expresses a geodesic flow on the group $\D^s(\T)$ for $s>3$ and any $b$. That, for $b=2$, this geodesic flow is compatible with the right-invariant two-dimensional $H^1$ metric is recalled in Appendix~\ref{app_geom1}. However, such a correspondence with a Riemannian structure works \emph{only} if $b=2$. This subtle result is given rigorous evidence in Theorem~\ref{thm_b=2} and we will now perform the proof. In accordance with the results of \cite{ES10,K09} we are working with the $b$-equation on the Fr\'echet-Lie group of smooth orientation-preserving diffeomorphisms of the torus.\\[.25cm]
\emph{Proof of Theorem~\ref{thm_b=2}}. We assume that, for a given $b\geq 2$ and $A\in\mathcal L_{\text{is}}^{\text{sym}}(C^{\infty}(\T))$, the
2D $b$-equation is the Euler equation on the torus
diffeomorphism group with respect to $\rho_\A$, i.e.,
$$u_t=-\A^{-1}\{u\cdot\nabla(\A u)+(\nabla u)^T\A u+\A u(\nabla\cdot u)\}.$$
The 2D $b$-equation reads
$$(\mathbb Lu)_t=-u\cdot\nabla (\mathbb Lu)-(\nabla u)^T(\mathbb Lu)-(b-1)(\mathbb Lu)(\nabla\cdot u).$$
Using that $(\mathbb Lu)_t=\mathbb Lu_t$ and resolving both equations with respect
to $u_t$ we get that
\bea\label{gl1}\A^{-1}\left\{u\cdot\nabla(\A u)+(\nabla u)^T\A u+\A u(\nabla\cdot u)\right\}&=&\nonumber\\
&&\hspace{-3cm}\mathbb L^{-1}\left\{u\cdot\nabla (\mathbb Lu)+(\nabla u)^T(\mathbb Lu)+(b-1)(\mathbb Lu)(\nabla\cdot u)\right\}.\eea
We now conclude that $\text{span}\{\mathbf 1\}$, $\mathbf 1=e_1+e_2$, is an invariant subspace of $\A$. We let $A1=\lambda$ and have $\A\mathbf 1=\lambda\mathbf 1$. Evaluating \autoref{gl1} for $e_1$ and $e_2$ shows that $\nabla\lambda=0$. The structure of \eqref{gl1} suggests that we have the freedom to scale the operator $\A$ so that we can assume that $\mathbf 1$ is a fixed point for $\A$. We next replace $u$ by $u+\sigma\mathbf 1$ in \eqref{gl1} and divide \eqref{gl1} by $\sigma$, to obtain for $\sigma\to\infty$ that
\beq\label{gl2}\A^{-1}\left[(\nabla u)^T+\nabla(\A u)+(\nabla\cdot u)\right]\mathbf 1=\mathbb L^{-1}\left[(\nabla u)^T+\nabla(\mathbb L u)+(b-1)(\nabla\cdot u)\right]\mathbf 1.\eeq
Let $n=(n_1,n_2)\in(2\pi\Z)^2\backslash\{(0,0)\}$, and write $z=(x,y)$ for the variable on $\T$. We will consider the functions $u_n=e^{\i nz}\mathbf 1$ in the following for which we have $\mathbb Lu_n=(1+n^2)u_n$ and $\mathbb L^{-1}u_n=(1+n^2)^{-1}u_n$, $n^2=n_1^2+n_2^2$. Set $v_n=\A u_n$. An explicit calculation of the left- and right-hand sides of \eqref{gl2} shows that we have the identity
\beq\label{gl3}\nabla v_n\cdot\mathbf 1-\i\alpha_{n}v_n=-\i\beta_n u_n\eeq
where
\begin{align}
\alpha_n=&\text{ diag}\left(\frac{n_1(b+1)+(b-1)n_2}{1+n^2}+n_1+n_2,\frac{n_2(b+1)+(b-1)n_1}{1+n^2}+n_1+n_2\right),\nonumber\\
\beta_n=&\text{ diag}\left(3n_1+n_2,3n_2+n_1\right).\nonumber\end{align}
Assume that $n_1\neq n_2$ first. Since $v_n=Ae^{\i nz}\mathbf 1$ we get that $Ae^{\i nz}=(1+n^2)e^{\i nz}$. For $n_1=n_2$ we see that the function $v_n=\frac{2}{b}(1+n^2)e^{\i nz}\mathbf 1$ constitutes a solution to \eqref{gl3}. We insert $u_n$ for $n_1=n_2\neq 0$ into \autoref{gl1} to get that $b=2$. Since $\{e^{\i nz};\,n\in(2\pi\Z)^2\}$ is a basis for $C^{\infty}(\T)$, it follows that $\A=\mathbb L$. This completes the proof of the theorem.
\hfill$\square$\\[.25cm]
We immediately obtain the following
\cor The geodesic flow for the two-dimensional Degasperis-Procesi equation
$$\frac{\partial m_i}{\partial t}=-u^j\frac{\partial m_i}{\partial x^j}-m_j\frac{\partial u^j}{\partial x^i}-2m_i\frac{\partial u^j}{\partial x^j},\quad m_i=(1-\Delta)u_i, \quad i=1,2,$$
on the diffeomorphism group of the torus is not related to any right-invariant metric on $\Diff^{\infty}(\T)$ with inertia operator $\text{\rm diag}(A,A)$, $A\in\L_{\text{\rm is}}^{\text{\rm sym}}(C^{\infty}(\T))$.
\endcor\rm
\section{Outlook}\label{sec_outlook}
It is natural to ask for $n$-dimensional generalizations of the $b$-equation which can be obtained by considering higher-order tensor densities and their Lie derivatives respectively. The results presented in Section~\ref{sec_2DCH} may easily extend to the $n$-dimensional case, but it is an open problem whether the geodesic flow on the $n$-torus is weakly Riemannian if and only if $b=2$.

Letting $\A=\text{diag}(\mu-\Delta,\mu-\Delta)$, where
$$\mu(u)=\int_\T u\d(x,y),$$
the family \eqref{2Db} becomes the two-dimensional $\mu$-$b$-equation. The special cases $b=2,3$, for which one obtains the two-dimensional $\mu$-Camassa-Holm and $\mu$-Degasperis-Procesi equations are mentioned in an appendix of the paper \cite{LMT10}. Since there is a significant interest in the one-dimensional partially averaged $b$-equation (see, e.g., \cite{EKK11,KLM08,K11,LMT10}), it is appealing to study geometric properties of the 2D $\mu$-$b$-equation. This will be part of a companion paper.

Since we restricted our attention to periodic functions in the present paper, we can now ask for a realization of the $b$-equation on diffeomorphism groups of non-compact manifolds, e.g., on the real axis. However, the group $\Diff^{\infty}(\R)$ of smooth and orientation-preserving diffeomorphisms $\R\to\R$ fails to be regular in the sense that not any element of its Lie algebra can be integrated into a one-parameter subgroup \cite{Mi06}. The geometric approach to non-periodic equations of type \eqref{2Db} is an open and attractive problem in the geometric analysis community where we probably need some deep and new ideas.
\begin{appendix}
\section{The geometric picture for the 2D Camassa-Holm equation}\label{app_geom1}
In this appendix, we provide some supplementary material concerning the geometric aspects of the 2D-CH equation. First, we give a full proof of Proposition~\ref{prop_geometry}.\\[.25cm]
\emph{Proof of Proposition~\ref{prop_geometry}.} That $\bar\nabla$ satisfies the properties (i)--(iv) follows immediately from our definitions. To prove smoothness, we verify that the map $\D^s(\T)\to\L^2_{\text{sym}}(H^s(\T);H^s(\T))$, $\phi\mapsto\Gamma_\phi$ is smooth. For $U\in H^s(\T)$ we write
$$G(\phi;U)=(\phi,\Gamma_\phi(U,U))=(\phi,(\A^{-1}P(U\circ\phi^{-1}))\circ\phi),$$
where $P$ is a polynomial differential operator. Introducing the maps
$$\tilde\A(\phi,V)=(\phi,R_\phi\circ\A\circ R_{\phi^{-1}}V),\quad\tilde P(\phi,V)=(\phi,R_\phi\circ P\circ R_{\phi^{-1}}V)$$
we see that $G=\tilde\A^{-1}\circ\tilde P$. Moreover $\tilde\A$ and $\tilde P$ are smooth maps, since $\A$ and $P$ are polynomial differential operators and $H^s(\T;\R)$ is a Banach algebra for $s>1$, \cite{A75}. It follows from
$$D\tilde\A_{(\phi,V)}=\left(
\begin{array}{cc}
\id & 0 \\
* & R_\phi\circ\A\circ R_{\phi^{-1}} \\
\end{array}
\right)
$$
that $D\tilde\A_{(\phi,V)}\colon H^s(\T)\times H^s(\T)\to H^s(\T)\times H^{s-2}(\T)$ is a bijective bounded linear map for any $(\phi,V)$. In view of the open mapping theorem and the inverse mapping theorem, $\tilde\A^{-1}$ is a smooth map and hence also $G$ is smooth.

We now turn our attention to the family $\ska{\cdot}{\cdot}_\phi$ of inner products. To prove that $(\D^s(\T),\ska{\cdot}{\cdot})$ is a Riemannian manifold, it remains to check that the map
$$\phi\mapsto Q(\phi;U,V)=\ska{U}{V}_\phi$$
is smooth for any $U,V\in H^s(\T)$. In view of the change of variables $(\tilde x,\tilde y)=\phi^{-1}(x,y)$, this follows immediately from the representation
$$Q(\phi;U,V)=\sum_{i=1}^2\int_\T\left\{U_iV_i+[\nabla U_i^T\cdot(\nabla\phi)^{-1}]\cdot[\nabla V_i^T\cdot(\nabla\phi)^{-1}]\right\}|\nabla\phi|\d(\tilde x,\tilde y).$$
Finally, to prove the compatibility of $\bar\nabla$ and $\ska{\cdot}{\cdot}$, we first prepare the identity
\beq\left.\frac{\dw}{\dw\eps}\right|_{\eps=0}Y_i(\phi+\eps X(\phi))\circ(\phi+\eps X(\phi))^{-1}=[DY_i(\phi)\cdot X(\phi)]\circ\phi^{-1}-\nabla v_i\cdot u,\label{auxformula}\eeq
where $X,Y\in\mathfrak X(\D^s(\T))$ and $u=X(\phi)\circ\phi^{-1}$, $v=Y(\phi)\circ\phi^{-1}$. Differentiating the relation $(\phi+\eps X(\phi))\circ(\phi+\eps X(\phi))^{-1}=\id$ with respect to $\eps$ at $\eps=0$, we get $\left.\frac{\dw}{\dw\eps}\right|_{\eps=0}(\phi+\eps X(\phi))^{-1}=-[(\nabla\phi)^{-1}X(\phi)]\circ\phi^{-1}$. Together with $\nabla(\phi^{-1})=(\nabla\phi)^{-1}\circ\phi^{-1}$, we obtain \eqref{auxformula}. Let $Z\in\mathfrak X(\D^s(\T))$ and write $w=Z(\phi)\circ\phi^{-1}$. In view of relation \eqref{auxformula}, we get that
\bea(X\ska{Y}{Z})(\phi)&=&\int_\T\left\{[DY(\phi)\cdot X(\phi)]\circ\phi^{-1}-\nabla v\cdot u\right\}\A w\d(x,y)\nonumber\\
&&+\int_\T\left\{[DZ(\phi)\cdot X(\phi)]\circ\phi^{-1}-\nabla w\cdot u\right\}\A v\d(x,y).\nonumber
\eea
On the other hand
$$\ska{\nabla_XY}{Z}_\phi=\int_\T\left\{[DY(\phi)\cdot X(\phi)]\circ\phi^{-1}-\Gamma(u,v)\right\}\A w\d(x,y)$$
so that it remains to prove the identity
\beq\int_\T[(\nabla v\cdot u)\A w+(\nabla w\cdot u)\A v]\d(x,y)=\int_\T[\Gamma(u,v)\A w+\Gamma(u,w)\A v]\d(x,y).\label{toprove}\eeq
In view of the definition of $\Gamma$ and the fact that the adjoint of $\nabla u$ with respect to the $L_2$ inner product is $(\nabla u)^T$, the right hand side of \autoref{toprove} equals
\bea &&-\frac{1}{2}\int_\T\{[(u\cdot\nabla)\A v]\cdot w+[(v\cdot\nabla)\A u]\cdot w +[(\nabla v)^T\A u]\cdot w+(\nabla\cdot v)(\A u)\cdot w\nonumber\\
&&\quad +(\nabla\cdot u)(\A v)\cdot w-[\A(\nabla v\cdot u)]\cdot w+(u\cdot\nabla)(\A w)\cdot v+(w\cdot\nabla)(\A u)\cdot v\nonumber\\
&&\quad +[(\nabla w)^T\A u]\cdot v+(\nabla\cdot w)(\A u)\cdot v+(\nabla\cdot u)(\A w)\cdot v-[\A(\nabla w\cdot u)]\cdot v\}\d(x,y).
\label{terms1}\eea
Using the identity
$$\int_\T\{(u\cdot\nabla)\A v+(\A v)(\nabla\cdot u)\}w\d(x,y)=-\int_\T\A v(u\cdot\nabla)w\d(x,y)$$
and $(u\cdot\nabla)v=\nabla v\cdot u$ we can rewrite \eqref{terms1} as
\bea&&\frac{1}{2}\int_\T\{2(\nabla v\cdot u)\A w+2(\nabla w\cdot u)\A v+(\nabla w\cdot v)\A u-(\nabla v)^T\A u\cdot w\nonumber\\
&&\quad +(\nabla v\cdot w)\A u-(\nabla w)^T\A u\cdot v\}\d(x,y).
\label{terms2}\eea
%
Simplifying these terms it follows that the expression \eqref{terms2} reduces to
$$\int_\T\{(\nabla v\cdot u)\A w+(\nabla w\cdot u)\A v\}\d(x,y)$$
and this is the left hand side of \autoref{toprove}. Altogether, this completes the proof of the proposition.
\hfill$\square$\\[.25cm]
Note that, as on general Banach manifolds, the geodesic flow $\phi(t;u_0)$ defines an exponential map $\exp_{\phi(0)}(u_0)=\phi(1;u_0)$.
We now explain briefly that the geodesic flow is indeed length minimizing in the sense that it minimizes the functional
$$L(\gamma)=\int_J\ska{\gamma_t(t)}{\gamma_t(t)}_{\gamma(t)}^{1/2}\d t.$$
In fact, the arguments are just a repetition of the general approach for Riemannian metrics on Banach manifolds, cf.~\cite{L95}, and they have been worked out by Lenells \cite{L07} for the 1D Camassa-Holm equation. Exchanging $\S$ with $\T$, we can proceed as in \cite{L07} to derive the following theorems.
\thm[Existence of normal neighborhoods] Let $\phi_0\in\D^s(\T)$. Given an open neighborhood $\mathcal V=\mathcal U_0\times B_\eps(0)$ of $(\phi_0,0)\in T\D^s(\T)$, there is an open neighborhood $\mathcal W\subset\mathcal U_0$ of $\phi_0$ in $\D^s(\T)$ such that any two points $\phi,\psi\in\mathcal W$ can by joined by a unique geodesic lying in $\mathcal U_0$, and such that for any $\phi\in\mathcal W$, the exponential map $\exp_\phi$ maps the open set in $T_\phi\D^s(\T)$ represented by $(\phi,B_\eps(0))$ diffeomorphically onto an open set $\mathcal U(\phi)$ containing $\mathcal W$.
\endthm\rm
The following theorem is proved by applying the Gauss lemma for the right-invariant metric $\ska{\cdot}{\cdot}_\phi$ on $\D^s(\T)$.
\thm Let $(\mathcal V,\mathcal W)$, $\mathcal V=\mathcal U_0\times B_\eps(0)$ constitute a normal neighborhood of an element $\phi_0\in\D^s(\T)$. Let $\alpha\colon[0,1]\to\D^s(\T)$ be the unique geodesic joining two points $\phi,\psi\in\mathcal W$. Then, for any piecewise $C^1$-path $\gamma\colon[0,1]\to\D^s(\T)$ joining $\phi$ and $\psi$, it holds that
$$L(\alpha)\leq L(\gamma).$$
If equality holds, then a reparametrization of $\gamma$ is equal to $\alpha$.\endthm\rm
Conversely, it holds globally that any length-minimizing curve is a geodesic.
\thm If $\alpha\colon[0,1]\to\D^s(\T)$ is a piecewise $C^1$-path parametrized by arc-length such that $L(\alpha)\leq L(\gamma)$ for all paths $\gamma$ in $\D^s(\T)$ joining $\alpha(0)$ and $\alpha(1)$, then $\alpha$ is a geodesic.
\endthm\rm
\section{The 2D Camassa-Holm equation as an Euler equation}\label{app_geom2}
In this appendix, we compare the geometric formalism for the 2D-CH to the geometric picture of the rigid body and the one-dimensional CH equation, cf.~also \cite{EKL11} where a similar comparison is worked out.

For a rigid body in $\R^3$ the configuration manifold for rotations around the origin is the finite-dimensional Lie group $G=SO(3)$ whose Lie algebra $\g=\mathfrak{so}(3)$ consists of antisymmetric $3\times 3$-matrices. The space $\mathfrak{so}(3)$ is canonically identified with $\R^3$ via the map
$$\hat{}\ :\R^3\to\mathfrak{so}(3),\quad x=(x_1,x_2,x_3)\mapsto\hat x=\left(%
\begin{array}{ccc}
  0 & -x_3 & x_2 \\
  x_3 & 0 & -x_1 \\
  -x_2 & x_1 & 0 \\
\end{array}%
\right).$$
We define a map $\check{}\,\colon\R^3\to\mathfrak{so}(3)^*$ by $(\check y,\hat x)=y\cdot x=y_1x_1+y_2x_2+y_3x_3$ to see that $\mathfrak{so}(3)^*\simeq\R^3$. Let $I:\mathfrak{so}(3) \to \mathfrak{so}(3)^*$ be the
inertia matrix of the body. A \emph{left}-invariant metric
$\ska{\cdot}{\cdot}$ on $SO(3)$ is defined by setting
$$\ska{a}{b}=a\cdot Ib, \qquad  a,b\in \R^3 \simeq \mathfrak{so}(3),$$
at the identity, and extending it to all of $SO(3)$ by left
invariance. The basic observation is that $R(t)$ is a geodesic on
$(SO(3),\ska{\cdot}{\cdot})$ if and only if
$\hat\Omega(t)=R(t)^{-1}\dot R(t)$ solves the classical
Euler equation for the motion of a rotating rigid body,
\beq\label{Eulerbody}
I\dot{\Omega} = (I\Omega)\times\Omega.
\eeq
Physically, $\hat\Omega(t)$ represents the angular velocity in a
frame of reference fixed with respect to the body. The angular
velocity in the spatially fixed frame is given by $\dot
R(t)R(t)^{-1}$. In other words: Applying left and right
translations to the material angular velocity $\dot R(t)$, one
obtains the \emph{body} and the \emph{spatial} angular velocities,
which are both elements of the Lie algebra $\mathfrak{so}(3)$. The
body and spatial angular momenta, which are elements of the dual
$\mathfrak{so}(3)^*$, are given by $\Pi(t)=I\Omega(t)$ and
$\pi(t)=R(t)\Pi(t)$, respectively. The body and spatial quantities
are related by  the adjoint and coadjoint actions
\begin{equation}\label{rigidbodyadjoints}
  \hat\omega(t)= \Ad_{R(t)}\hat\Omega(t)=R(t)\hat\Omega(t)R(t)^{-1},\qquad\Pi(t)=\Ad^*_{R(t)}\pi(t).
\end{equation}
Conservation of (spatial) angular momentum implies that $\pi$ is in fact constant in
time, i.e.
\begin{equation}\label{piconserved}
  \frac{\dw\pi}{\dw t} = 0.
\end{equation}
Finally, the Euler equation \eqref{Eulerbody} can be rewritten in the form
$$\check\Pi_t=\ad_{\hat{\Omega}}^*\check\Pi.$$
It is shown in \cite{K04} that the geometric approach works also well for the 1-dimensional periodic Camassa-Holm equation \eqref{CH}
for which the configuration space is $G = \Diff(\S)$ with multiplication
$(\varphi,\psi)\mapsto\varphi\circ\psi$. Elements of the Lie
algebra $\mathfrak{g}$ are identified with functions $\S \to \R$.
A \emph{right}-invariant metric is defined by setting
$$\ska{u}{v}_{H^1}=\int_{\S} uAv\dx=\int_{\S}(uv+u_xv_x)\dx,$$
where $A=1-\partial_x^2\colon\mathfrak{g} \to \mathfrak{g}^*$ is the
inertia operator and $\g^*\simeq\g$ is a regular part of the topological dual space $\g'$. The basic observation is that $\varphi(t)$ is a
geodesic in $(\Diff(\S),\ska{\cdot}{\cdot}_{H^1})$ if and only if
$u(t)=DR_{\varphi(t)^{-1}}\varphi_t(t)=\varphi_t(t)\circ\varphi(t)^{-1}$
satisfies (\ref{CH}). In other words, the CH equation is the
Euler equation on $(\Diff(\S),\ska{\cdot}{\cdot}_{H^1})$. Letting
$U = DL_{\varphi^{-1}}\varphi_t = (u\circ\varphi)\varphi_x^{-1}$,
$U$ and $u$ are the analogs of the body and spatial angular
velocities: they are obtained by left and right translation, respectively, of
the material velocity $\varphi_t$ to the Lie algebra. The momentum
in the spatial frame is $m=Au$. The analog of equation
(\ref{rigidbodyadjoints}) is
$$u(t) = \Ad_{\varphi(t)} U(t), \qquad m_0(t) = \Ad_{\varphi(t)}^*m(t),$$
where $m_0 = (m\circ\varphi)\varphi_x^2$ is the momentum in the
body frame. Since the metric now is right-invariant instead of
left-invariant, the analog of the conservation law
(\ref{piconserved}) is that the momentum $m_0$ in the body frame
is conserved,
$$\frac{\dw m_0}{\dw t} = 0, \qquad \text{i.e.} \qquad (m \circ \varphi)\varphi_x^2 = m_0.$$
The Camassa-Holm equation \eqref{CH} in terms of the Eulerian velocity $u$ reads
$$u_t=-\ad_u^*u.$$
In the present appendix, we show that the two-dimensional periodic Camassa-Holm equation can be described within Arnolds geometric approach.
For the 2D-CH equation (\ref{2DCH}) the configuration space is the
group $G = \Diff(\T)$ introduced in Section~\ref{sec_preliminaries}. The Lie algebra
$\mathfrak{g}$ is identified with the set of functions $\T \to \R^2$. The inertia operator is $\A=\text{diag}(1-\Delta,1-\Delta)$
and the metric is the right-invariant metric $\langle \cdot, \cdot
\rangle$ defined in (\ref{metric}). The basic observation is
that $\phi(t)$ is a geodesic in $(\Diff(\T),\ska{\cdot}{\cdot})$ if and only if
$$u = DR_{\varphi(t)^{-1}}\varphi_t(t)=\phi_t\circ\phi^{-1}$$
satisfies \autoref{2DCH}. The analog of the body angular velocity is $$U = DL_{\varphi^{-1}}\varphi_t=(\nabla\phi)^{-1}\phi_t.$$
The spatial momentum is given by $m = \A u$.
The analog of equation (\ref{rigidbodyadjoints}) reads
$$u(t) = \Ad_{\varphi(t)} U(t)$$
and
$$m_0= \Ad_{\varphi(t)}^*m(t)=(\nabla\phi)^T(m\circ\phi)|\nabla\phi|,$$
where $m_0$ is the momentum in the body frame. The analog of the conservation law (\ref{piconserved}) is that the momentum $m_0$ in the body frame is conserved. It is left to the reader to perform the explicit calculations showing that the derivative of $m_0$ is indeed zero.
\noindent Finally, in view of \eqref{ad} and \eqref{2DCHEuler}, we can recast the 2D Camassa-Holm equation in terms of the Eulerian variable $u$ as
$$u_t=-\ad_u^*u.$$
Let us summarize our results in the following tabular.
\begin{center}
\begin{scriptsize}
\begin{tabular}{|l|c|c|c|}
  \hline
   & Rigid body & CH & 2D-CH \\
  \hline
  configuration space &$SO(3)$&$\Diff(\S)$&$\Diff(\T)$\\
  material velocity & $\dot R$ & $\varphi_t$ & $\phi_t=(\varphi_{1t},\varphi_{2t})$\\
  spatial velocity & $\hat\omega=\dot RR^{-1}$ & $u = \varphi_t \circ \varphi^{-1}$ & $u = \phi_t\circ \varphi^{-1} $ \\
  body velocity & $\hat\Omega=R^{-1}\dot R$ & $U=\frac{\varphi_t}{\varphi_x}$ & $U=(\nabla\phi)^{-1}\phi_t$ \\
  inertia operator & $I$ & $A=1-\partial_x^2$ & $\A=\begin{pmatrix}
  1-\Delta & 0 \\
  0 & 1-\Delta \end{pmatrix}$ \\
  spatial momentum & $\pi=R\Pi$ & $m=Au$ & $m =\A u$ \\
  body momentum & $\Pi=I\Omega$ & $m_0=(m\circ\varphi)\varphi_x^2$ & $ m_0 = (\nabla\phi)^T(m\circ\phi)|\nabla\phi|$  \\
  spatial velocity (Ad) & $\hat\omega = \Ad_{R}\hat\Omega$ & $u = \Ad_{\varphi} U$ & $u = \Ad_{\varphi} U$\\
  body momentum (Ad*) & $\Pi=\Ad^*_{R}\pi$  & $m_0 =\Ad^*_\varphi m$ & $m_0 =\Ad_{\varphi}^*m$ \\
  momentum conservation & $\pi = \text{const.}$ & $m_0 = \text{const.}$ & $m_0 = \text{const.}$ \\
  \hline
\end{tabular}
\end{scriptsize}
\end{center}
\end{appendix}
\vspace{.5cm}

\end{document}